\documentclass[aps,preprint,superscriptaddress,floatfix,amsmath,amssymb]{revtex4}
\usepackage{inputenc}
\usepackage[a-1b]{pdfx}
\usepackage{amsmath}
\usepackage{amssymb}
\usepackage{amsthm}
\usepackage{algcompatible}
\usepackage{algpseudocode}
\usepackage{newfloat}
\usepackage{bm}
\usepackage{graphicx}
\usepackage{subfig}
\usepackage{hyperref}

\begin{document}

\title{Approximation of Optimal Control Surfaces for the Bass Model with Stochastic Dynamics}

\author{Gabriel Nicolosi}
\email{gxr286@psu.edu}
\affiliation{Harold and Inge Marcus Dept.\ of Industrial and Manufacturing Engineering, Penn State University,
University Park, PA 16802}

\author{Christopher Griffin}%
\email{griffinch@psu.edu}
\affiliation{Applied Research Laboratory, Penn State University,
University Park, PA 16802}

\date{April 19, 2023}

\begin{abstract}
{\small The Bass diffusion equation is a well-known and established modeling approach for describing new product adoption in a competitive market. This model also describes diffusion phenomena in various contexts: infectious disease spread modeling and estimation, rumor spread on social networks, prediction of renewable energy technology markets, among others. Most of these models, however, consider a deterministic trajectory of the associated state variable (\textit{e.g.}, market-share). In reality, the diffusion process is subject to noise, and a stochastic component must be added to the state dynamics. The stochastic Bass model has also been studied in many areas, such as energy markets and marketing. Exploring the stochastic version of the Bass diffusion model, we propose in this work an approximation of (stochastic) optimal control surfaces for a continuous-time problem arising from a $2\times2$ skew symmetric evolutionary game, providing the stochastic counter-part of the Fourier-based optimal control approximation already existent in the literature.}
\end{abstract}

\maketitle

\section*{Keywords}
Optimal control, Fourier approximation, evolutionary game, learning, stochastic Bass model.

\section{Introduction and Relevant Literature}
With its origins in the marketing literature, the Bass diffusion equation is a classical modeling approach initially developed to describe the underlying dynamics of new product adoption in a competitive market \cite{B69}. Its usefulness in modeling nonlinear phenomena is evident by its later application in diverse areas: infectious disease spread modeling and estimation \cite{EDDT21}, product sales forecasting \cite{FCC17}, renewable energy technology diffusion \cite{KA16}, among others. The Bass model has also been used in the realm of non-deterministic problems: electricity markets \cite{SG97} and sales forecasting \cite{GK19} are just two examples. In the context of digital product maintenance, \cite{FG17} derived a deterministic Bass-constrained optimal control problem (OCP) from a skew-symmetric game (a game with skew-symmetric payoff matrix), for which a solution was found by means of the classical Pontryagin's Maximum Principle (PMP), a classical result in the theory of optimal control that provides necessary conditions for control optimality. More recently, \cite{NFG22} proposed a backpropagation inspired Fourier-based approximation procedure for the same problem, approximating the solution found in \cite{FG17} and generalizing it over a predefined range of system initial conditions. The work in \cite{NFG22} falls under the wide umbrella of global orthogonal collocation in numerical methods for OCPs, a well-established and active field of research. See \cite{B20} for an overview. Using the same approach as in \cite{NFG22}, we construct approximating optimal control surfaces for the stochastic variation of the OCP arising from a $2\times2$ skew-symmetric evolutionary game, constrained by the Bass equation. In order to introduce stochasticity to the original problem, noise is added to the state dynamics in two different forms: additive noise to the fully observable state $x(t)$, and additive noise to the observer $y(t)$. This approach of adding noise is consistent with the literature of stochastic control and realistically models the underlying uncertainty involved in the diffusion process . Further details are discussed in Section~\ref{sec: Preliminaries and Problem Construction}, in which both deterministic and stochastic versions of the problem in question are introduced. Section~\ref{sec: OC Surface Approximation} introduces the Fourier-based control parameterization used to generate the approximating optimal control and state surfaces. Some results from computational experiments are shown in Section~\ref{sec: Computational Results}. Finally, the general context of this work and future directions are briefly discussed in Section~\ref{sec: Context and Future Directions}.  

\section{Preliminaries and Problem Construction}\label{sec: Preliminaries and Problem Construction}
Consider the evolutionary dynamics of a population consisting of $n_1 + n_2 = n$ individuals, wherein these two groups represent non-consumer ($n_1$) and consumer ($n_2$) partitions of the total population. Defining $x_i=\frac{n_i}{n}$ as the percentages of individuals in each group, the system of differential equations describing the evolution of each group and the underlying market-share dynamics is
\begin{equation}
    \begin{aligned}
            \dot{x}_1 &= -\rho x_1  x_2\\
            \dot{x}_2 &= \rho x_2 x_1
    \label{eqn: replicator dynamics}
    \end{aligned}
\end{equation}
where $\rho$ is the interaction payoff between both population partitions, giving rise to the skew-symmetric pay-off matrix 
$A=\bigl(\begin{smallmatrix} 
0    & -\rho \\
\rho & 0     \\
\end{smallmatrix}\bigl)$. 
We note that Equation~(\ref{eqn: replicator dynamics}) is an instance of the so-called replicator dynamics \cite{SS83}. Since $x_1 + x_2 = 1$, Equation~(\ref{eqn: replicator dynamics}) is reduced to $\dot{x} = \rho x(1-x), \label{eqn: Bass}$
which, in its turn, is a special case of the Bass diffusion equation from \cite{B69}, \textit{i.e.}:
\begin{equation}
    \dot{x} = p(1-x) + q(1-x)x
\end{equation}
for when the so-called coefficient of innovation $p=0$.
Now, consider an extraneous mechanism $u = u(t)$ that seeks to alter the dynamics shown in Equation~(\ref{eqn: Bass}), such as product maintenance effort aiming increased revenue stream through client base expansion. The Bass dynamics then becomes
\begin{equation}
    \dot{x} = \pi(u)x(1-x),
\end{equation}
where $\pi(u) = \beta u - \xi$ is a linear function of $u$, which is positive if and only if $\beta u > \xi$, \textit{i.e.}, the intervention effect on the market-share dynamics is only positive against the next best alternative $\xi$. Following the example in \cite{FG17}, a company wishes to minimize a continuous performance index from the launching of a digital product at $t=0$ up until the end of its life-cycle at $t=T$:
\begin{equation}
    \min_u \int_0^T (Cu^2 - \alpha ux) \; dt,
    \label{eqn: performance index}
\end{equation}
where $Cu^2$ is the quadratic cost of maintenance effort and $\alpha ux$ is the linear revenue stream. If we require the solution of Equation~(\ref{eqn: performance index}) to respect Equation~(\ref{eqn: Bass}) departing from a given $x(0)=x_0$, we have a OCP with a deterministic Bass dynamics constraint.  
\subsection{The Deterministic Bass Equation Constrained OCP}
The OCP of interest, with dynamics described by a deterministic Bass equation, is then formulated as:
\begin{equation}
    \begin{aligned}
            \min_u &\int_0^T (Cu^2 - \alpha u x) \; dt\\
            s.t. \quad
            &\dot{x} = x(1-x)(\beta u - \xi)\\
            &x(0) = x_0\\
            &u(t) \ge 0.
        \end{aligned}
\label{eqn: deterministic OCP}
\end{equation}
The solution to this OCP can be found analytically by applying the PMP like done in \cite{FG17}. If instead an optimal control solution surface is desired for a range of initial conditions $x_0 \in X_0$, the procedure introduced by \cite{NFG22} can be applied, yielding approximations based on a two-dimensional Fourier series, dependent both on time $t$ and initial condition $x_0$. The goal of this paper, however, is to provide an approximation in the taste of \cite{NFG22} to the stochastic variation of OCP~(\ref{eqn: deterministic OCP}), introduced next.

\subsection{The Stochastic Bass Equation Constrained OCP}
Building on OCP~(\ref{eqn: deterministic OCP}), we propose two stochastic model variations for the problem introduced above. The first model considers that the state dynamics is subject to additive white noise $\xi(t)$, giving a stochastic Bass model described by
\begin{equation}
    \dot{x} = x(1-x)(\beta u - \eta) + \sigma \xi(t),
    \label{eqn: StochasticBass model 1}
\end{equation}
where $\sigma \in \mathbb{R}$. The second variation assumes deterministic dynamics as in \cite{NFG22}, which in its turn cannot be observed without the presence of noise. Defining the observed state as $y(t)$, the noisy (observed) trajectory is formulated as
\begin{equation}
    y(t) = x(t) + \sigma\xi(t).
    \label{eqn: StochasticBass model 2}
\end{equation}
In what follows, we present the optimal control problem (OCP) formulation for both model dynamics variations presented in Equations~(\ref{eqn: StochasticBass model 1}) and~(\ref{eqn: StochasticBass model 2}). In light of the first stochastic model described above, the OCP under stochastic dynamics becomes
\begin{equation}
\begin{aligned}
    \min_{u} \quad & E\left[\int_0^T \left(Cu^2 - \alpha x u\right) \; dt \right]\\
            s.t. \quad &\dot{x} = x(1-x)(\beta u - \eta) + \sigma \xi(t)\\
                    &x(0) = x_0\\
                    &u \in L_2([0,T]),
\end{aligned}
\label{eqn: OCP1}                    
\end{equation}
where the expectation operator $E$ is needed in the objective functional in order to address the nondeterministic trajectory $x(t)$. The second OCP is defined in the same manner, however observing the trajectory, derived from Equation~(\ref{eqn: StochasticBass model 1}), through Equation~(\ref{eqn: StochasticBass model 2}). We also note that Equation~(\ref{eqn: StochasticBass model 1}) can be rewritten as
\begin{equation}
    dx = x(1-x)(\beta u - \eta)\;dt + \sigma\;dW,
\end{equation}
which is consistent with the notation in the literature of stochastic differential equations, and $dW$ comes from the fact that the Brownian motion $W(t)$ has its time derivative defined as a white noise process, that is, $\dot{W}(t) = \xi(t) \sim \mathcal{N}\left(0, \sigma\right)$.
Next, we introduced the control approximation $u \approx \hat{u}(t, x_0)$ and formulate the stochastic OCP, such as in (\ref{eqn: OCP1}), as a nonlinear programming (NLP) problem to be solved for the Fourier coefficients of the control approximation.  

\section{Optimal Control Surface Approximation}\label{sec: OC Surface Approximation}
Our goal is to find a control approximation $\hat{u}(t, x_0)$ that solves the OCP~(\ref{eqn: OCP1}) for any arbitrary initial condition in a predetermine set $X_0$. For this, we adapt the approximation scheme proposed in \cite{NFG22}. A two-variable cosine Fourier series is used to approximate the controllers $u(t, x_0)$:
\begin{equation}
    u(t, x_0) \approx \hat{u}(t, x_0) = \sum_{m=0}^M\sum_{n=0}^N a_{mn} \cos\left(\frac{m\pi t}{T}\right) \cos\left(\frac{n\pi x_0}{\mathcal{X}_0}\right)
    \label{eqn: Fourier}
\end{equation}
where $a_{mn}$ are Fourier coefficients, with $M, N \in \mathbb{N}$, $T$ is the same terminal time from~(\ref{eqn: OCP1}) and $\mathcal{X}_0$ is the range within which the initials conditions vary in $X_0$. The objective functional in~(\ref{eqn: OCP1}) is then transformed into a finite dimensional optimization problem, where the trajectory $x$ satisfies its original dynamics constraint, given an initial condition, wherein $u$ is approximated as in Equation~(\ref{eqn: Fourier}). 

\subsection{Stochastic Dynamics}\label{subsec: Stochastic Dynamics}
Substituting the approximation~(\ref{eqn: Fourier}) into~(\ref{eqn: StochasticBass model 1}), we obtain
\begin{equation}
    \dot{x} = x(1-x)(\beta \hat{u} - \eta) + \sigma \xi(t)
    \label{eqn: ODE}
\end{equation}
with $x(0)=x_0$, that is, the original dynamics constraint in OCP~(\ref{eqn: OCP1}) wherein $u$ is approximated as in Equation~(\ref{eqn: Fourier}). For the stochastic constrained problem addressed in this work, the non-deterministic component of~(\ref{eqn: ODE}) must be taken into consideration. Therefore, for a given initial condition, if $\varphi(t)$ satisfies Equation~(\ref{eqn: StochasticBass model 1}) without the presence of noise $\xi(t)$, then $\widetilde{\varphi}(t) = \varphi(t) + \sigma W(t)$ satisfies its stochastic counterpart, where $\sigma W(t)$ is the Brownian motion affecting the state trajectory. The new objective functional becomes:
\begin{equation}
      E\left[ \int_0^T \left(Cu^2 - \alpha\left(\varphi + \sigma W(t)\right)u\right) \; dt\right]. 
      \label{eqn: stoch. obj. functional}
\end{equation}
However, the new objective functional is equivalent to its deterministic counterpart. To show this, Equation~(\ref{eqn: stoch. obj. functional}) is rewritten as:
\begin{equation}
    \int_0^T \left(Cu^2 - \alpha\varphi u \right)\;dt  - \alpha\sigma u\;E\left[\int_0^T  W(t)  \;dt\right] = \int_0^T \left(Cu^2 - \alpha\varphi u \right)\;dt,
    \label{eqn: stoch. obj. functional2}
\end{equation}
since, on average, the definite integral of Brownian motion $W(t)$ is zero. Therefore, in this case, the objective functional remains unchanged from the deterministic one. Thus, the optimization procedure can be used exactly as in \cite{NFG22}. 

\subsection{Stochastic State Observer}\label{subsec: Stochastic State Observer}
Now we address the case in which the state trajectory is observed as described by Equation~(\ref{eqn: StochasticBass model 2}). For this, we substitute the right-hand side of Equation~(\ref{eqn: StochasticBass model 2}) into our objective functional, obtaining:
\begin{equation}
    \begin{aligned}
        E\left[\int_0^T \left(Cu^2 - \alpha \left(\varphi + \sigma \xi(t)\right) u\right)\;dt\right]
        = \int_0^T \left(Cu^2 - \alpha\varphi u \right)\;dt - \alpha\sigma u\;E\left[\int_0^T  \xi(t)  \;dt\right].
    \end{aligned}
\end{equation}
By definition, $\dot{W}(t) = \xi(t) $, therefore 
\begin{equation}
    \int_0^T  \xi(t)  \;dt = \int_0^T  \dot{W}(t)  \;dt = W(T) - W(0) = W(T).
\end{equation}
Thus, since $E[W(T)] = 0$, once again, the objective functional reduces to the deterministic case. We conclude that both in this case, and as in~(\ref{subsec: Stochastic Dynamics}), the control approximation~(\ref{eqn: Fourier}) is robust to the types of noise considered therein, providing the same approximation capabilities previously demonstrated in \cite{NFG22}.

\subsection{Nonlinear Programming Problem}
Using the control approximation from Equation~(\ref{eqn: Fourier}), we construct the following problem:
\begin{equation}
    \begin{aligned}
        \min_{u} \quad &E\left[\int_0^T \left(Cu^2 - \alpha x u\right) \; dt \right]\\
        s.t.\;\; &\dot{x} = x(1-x)(\beta \hat{u} - \eta) + \sigma \xi(t)\\
        &x(0) = x_0 \\
        &u(t,x_0) \approx \sum_{m=0}^M \sum_{n=0}^N a_{mn} \cos\left(\frac{m\pi t}{T}\right) \cos\left(\frac{n\pi x_0}{L}\right).
    \end{aligned}
    \label{eqn: NLP}
\end{equation}
We note that a similar formulation is valid for the case considered in Subsection~(\ref{subsec: Stochastic State Observer}). Like discussed in Subsection~(\ref{subsec: Stochastic Dynamics}), if $\varphi$ satisfies the flow determined by the deterministic component of the Bass dynamical constraint for a given initial condition $x_0$, then the OCP~(\ref{eqn: NLP}) becomes the following unconstrained NLP:
\begin{equation}
     \min_\mathbf{a} \quad J(\mathbf{a};\;x_0) = \int_0^T \left(C\hat{u}(\mathbf{a})^2 - \alpha\varphi(\mathbf{a};\;x_0)\hat{u}(\mathbf{a}) \right)\;dt.
     \label{eqn: NLP2}
\end{equation}
Finally, the approximating surface of optimal controllers is found by solving
\begin{equation}
    \min_{\mathbf{a}} \quad \sum_{x_0 \in X_0} J(\mathbf{a};\;x_0),
    \label{eqn: NLP3}
\end{equation}
which can be done by a first order method like gradient descent, taking place in the space of Fourier coefficients $\mathbf{a}$.
\section{Computational Results}\label{sec: Computational Results}
\begin{figure}[htpb]
     \centering
     \includegraphics[width=0.47\textwidth]{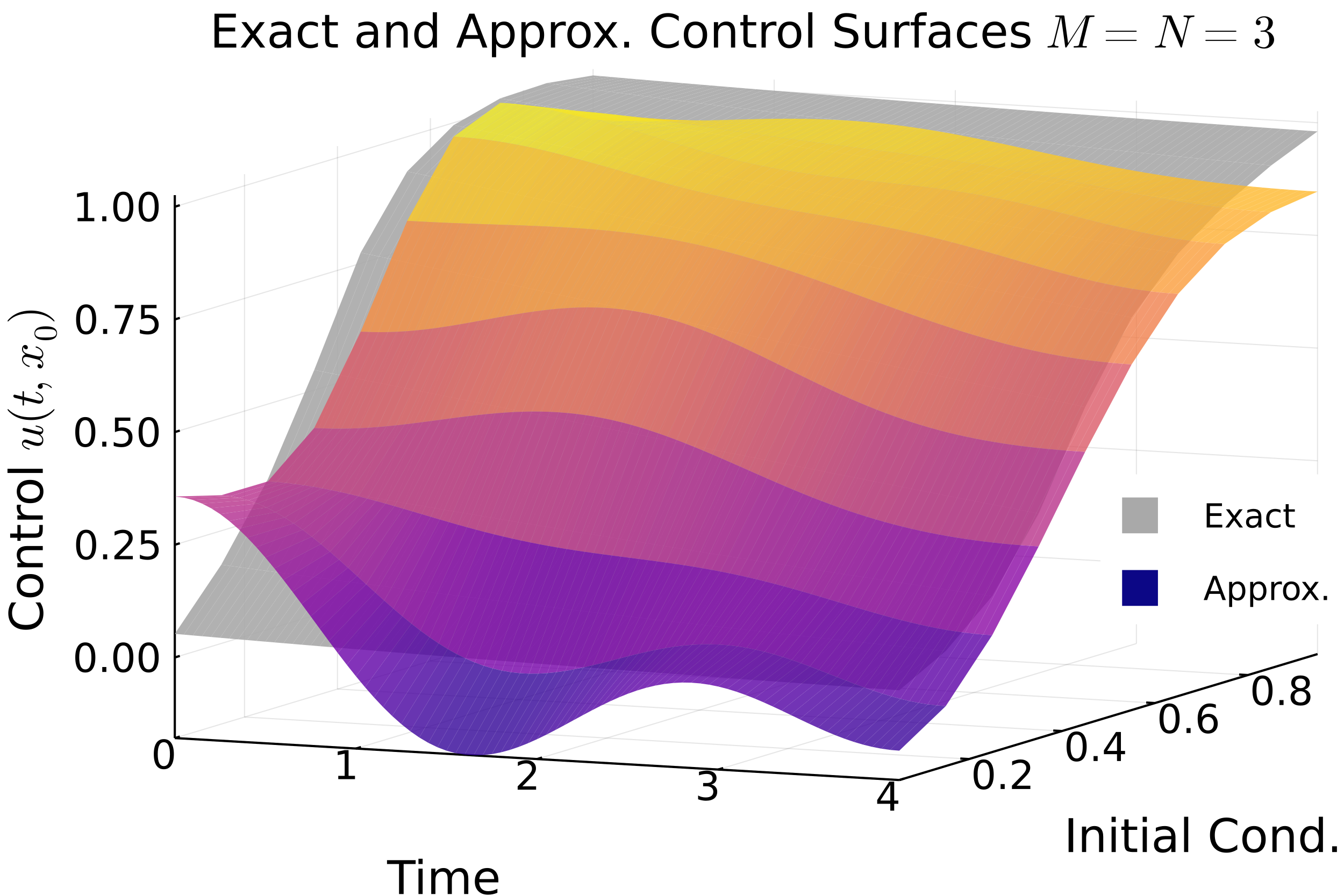}
     \qquad
     \includegraphics[width=0.47\textwidth]{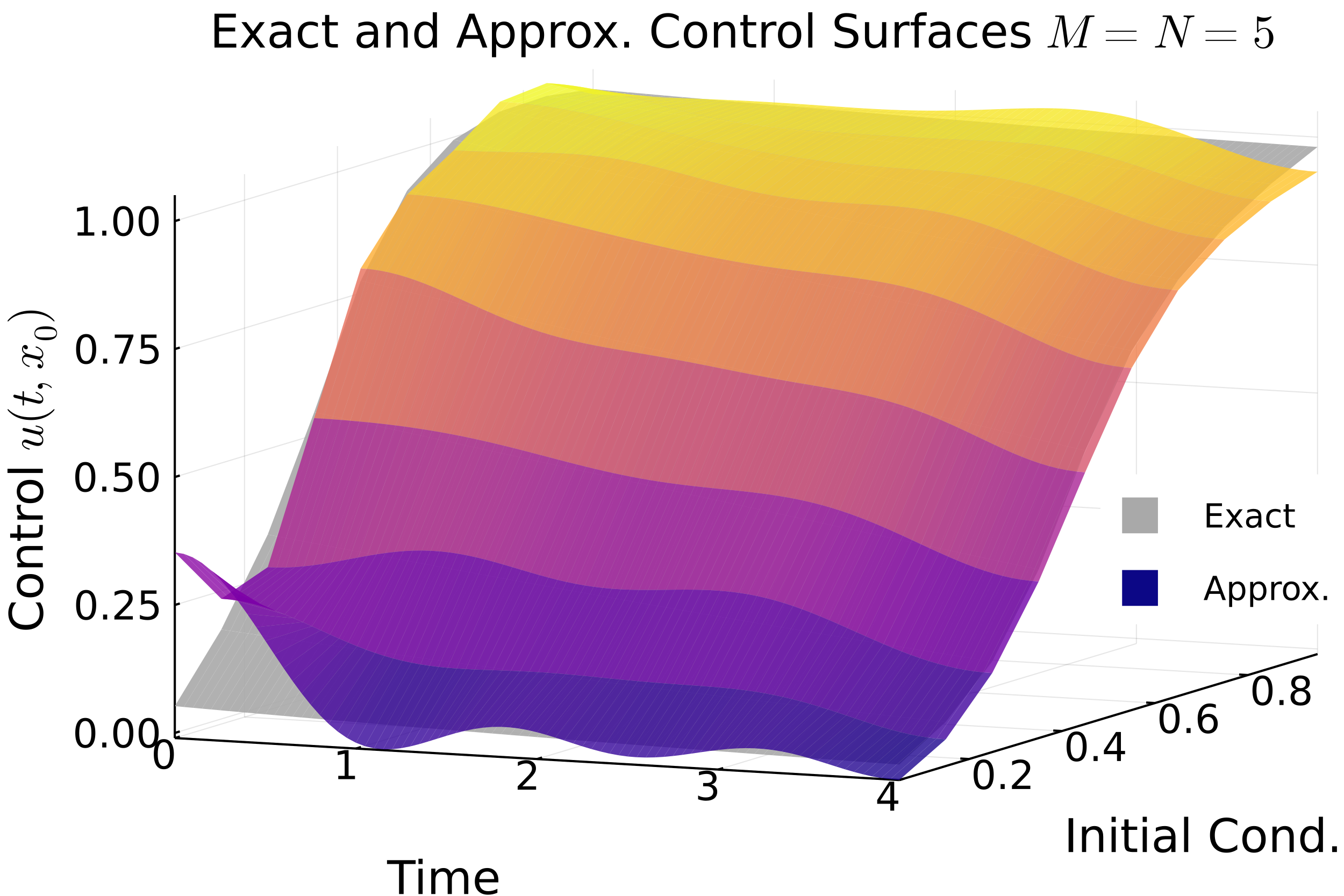}
     \qquad
     \includegraphics[width=0.47\textwidth]{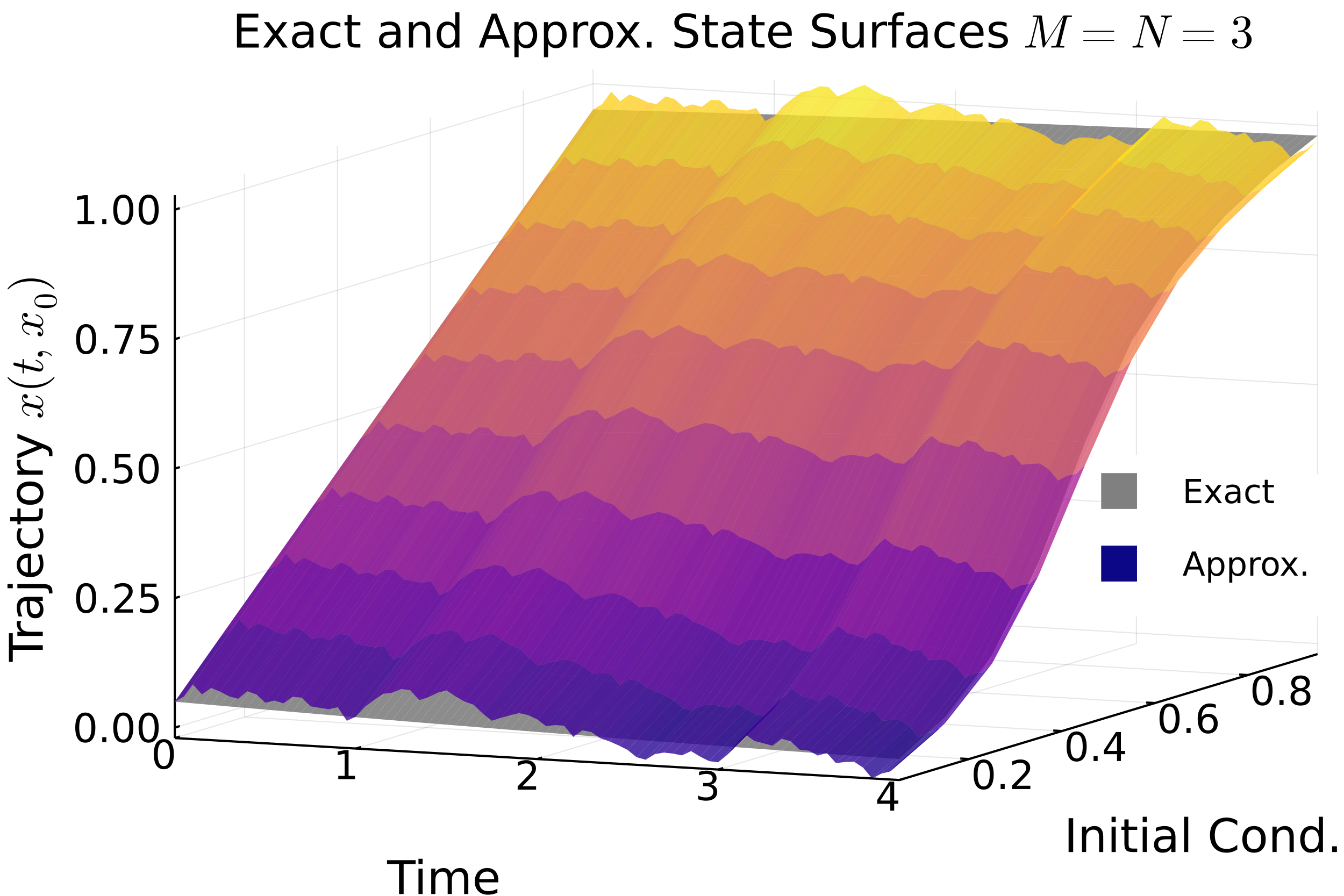}
     \qquad
     \includegraphics[width=0.47\textwidth]{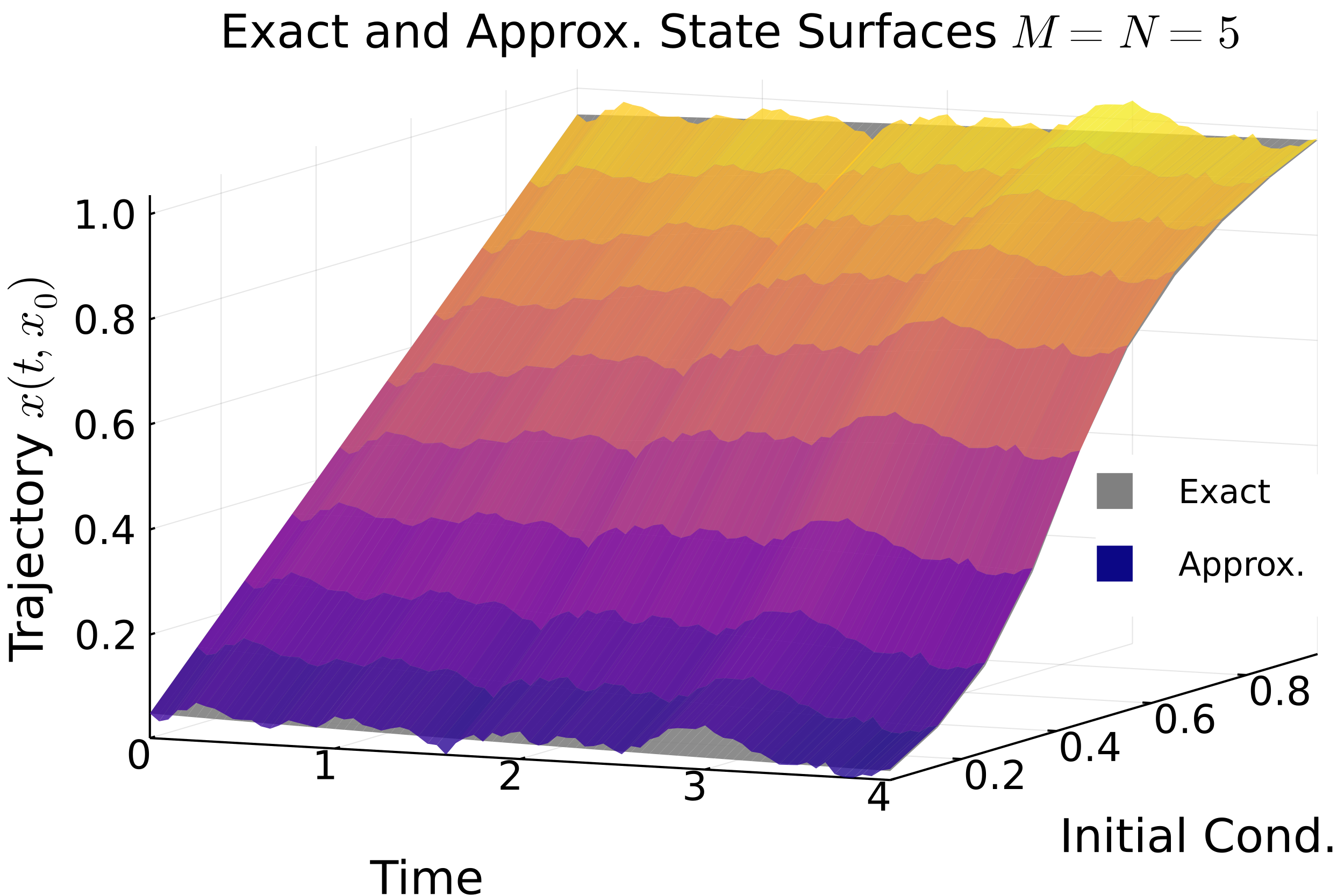}         
     \qquad
     \includegraphics[width=0.47\textwidth]{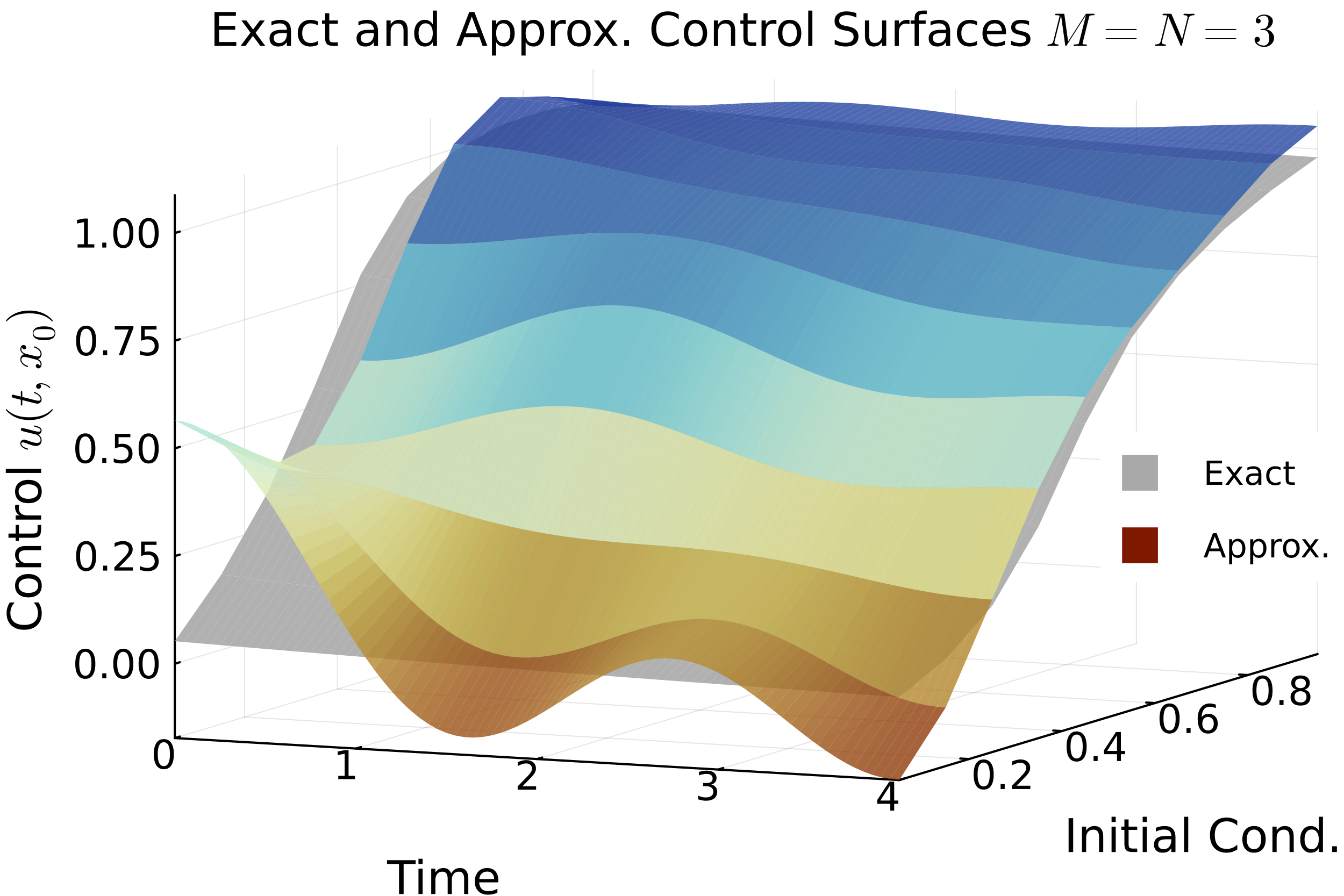}
     \qquad
     \includegraphics[width=0.47\textwidth]{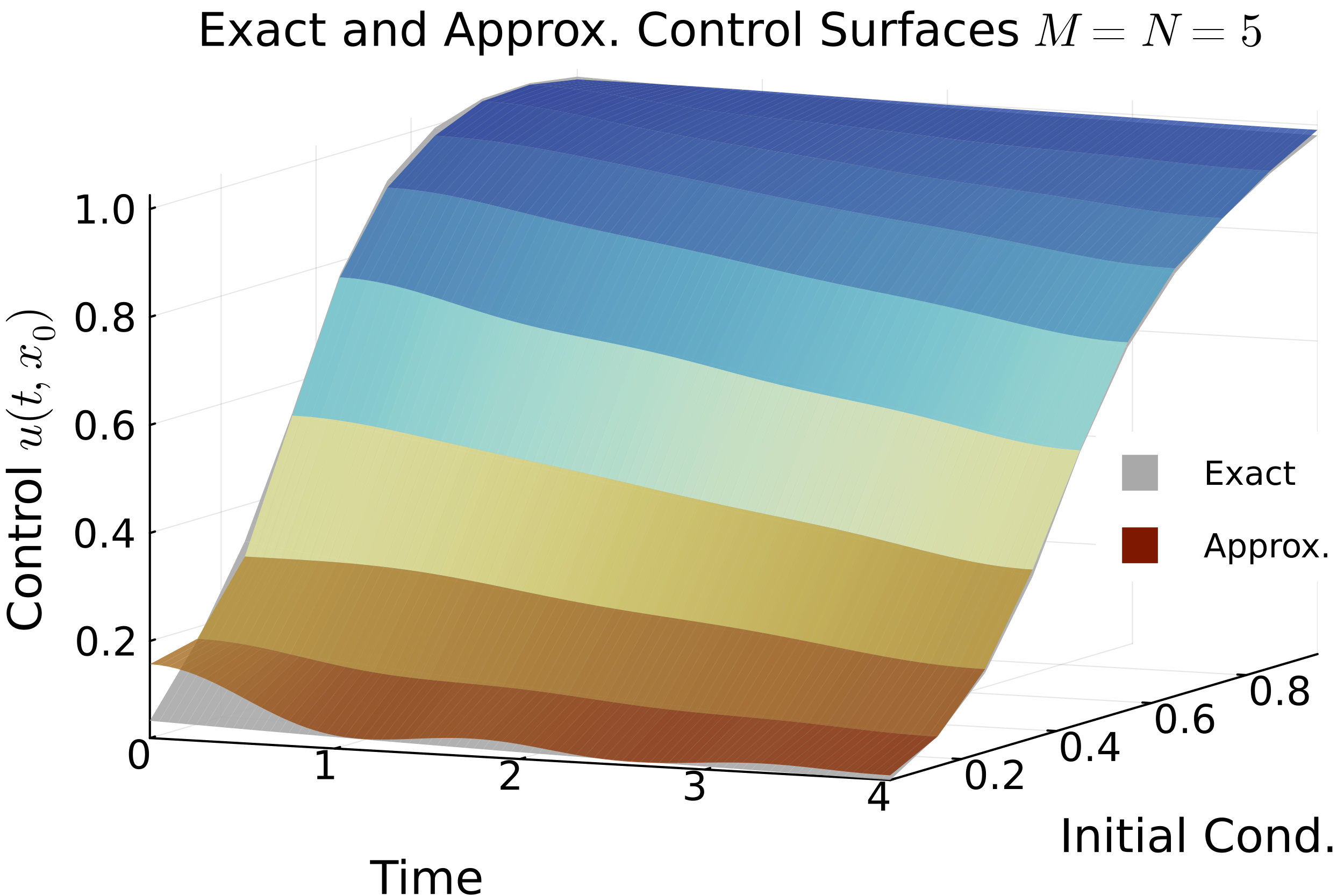}         
     \qquad
     \includegraphics[width=0.47\textwidth]{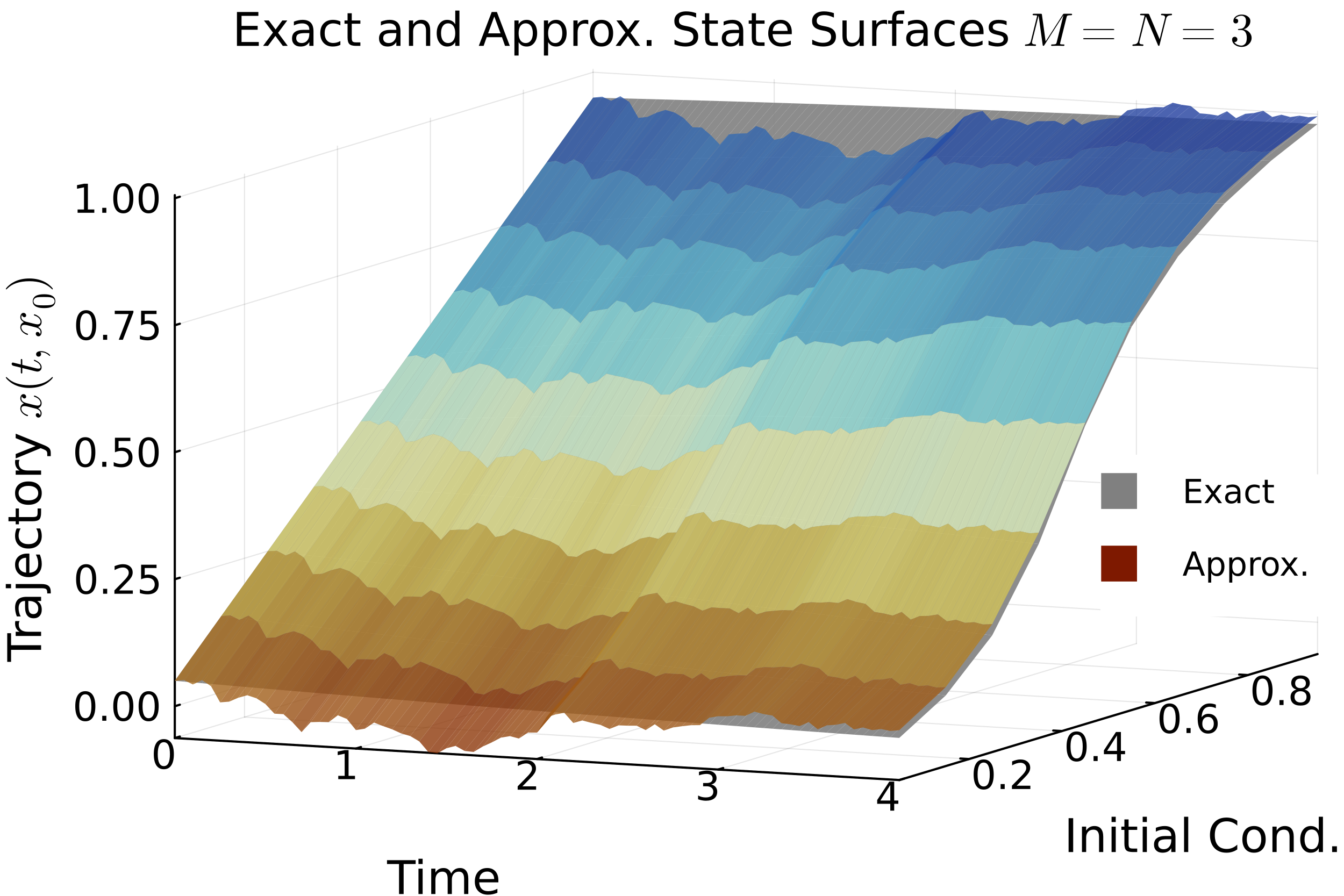}
     \qquad     
     \includegraphics[width=0.47\textwidth]{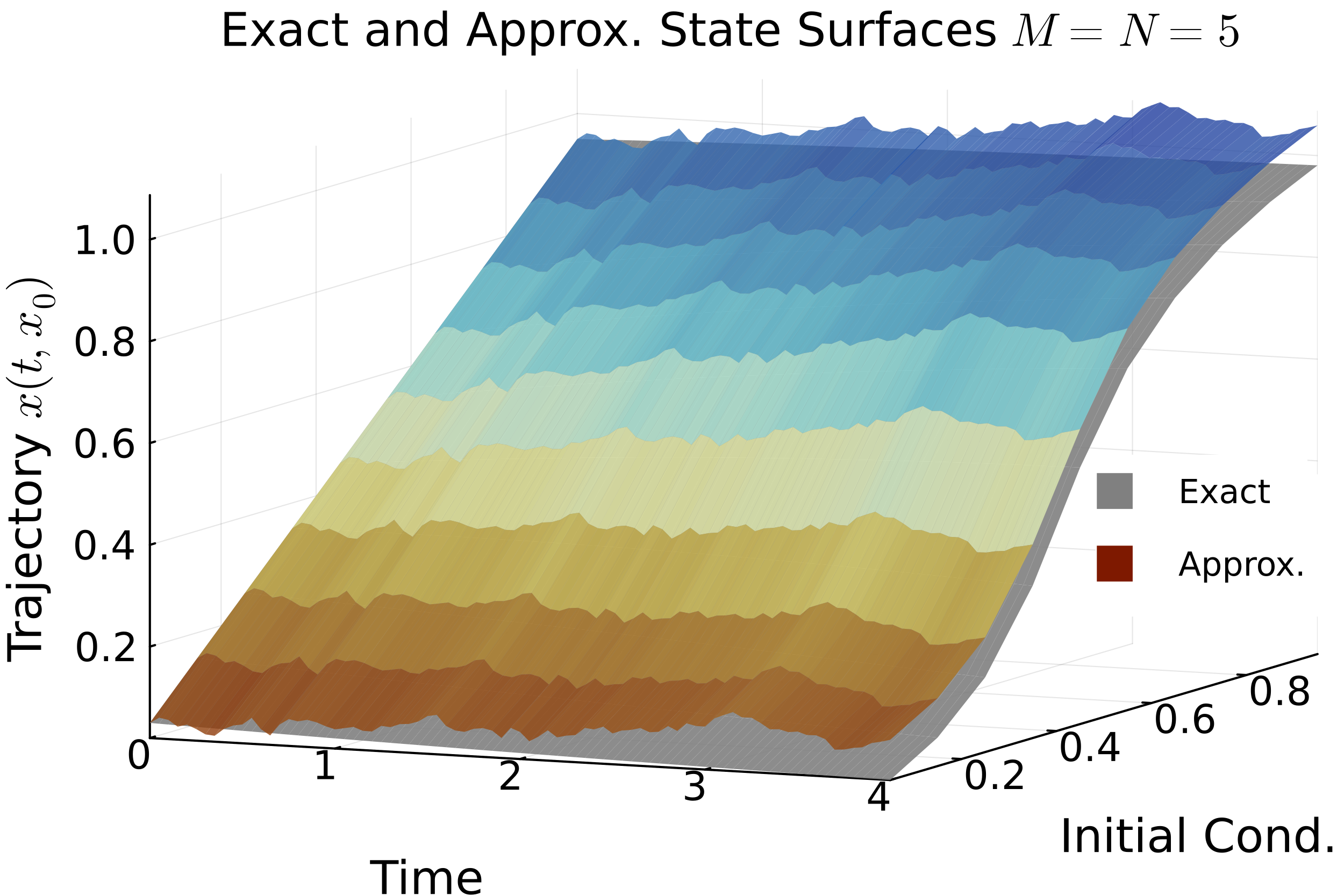}
\caption{Optimal control surface approximations and (sample) optimal state trajectory surfaces for two choices of Fourier coefficients: $M=N=3$ and $M=N=5$ for the case of stochastic dynamics $\left((a)-(d)\right)$ and noisy observer $\left((e)-(h)\right)$ for $\sigma=0.1$.}
\label{fig: Surfaces}   
\end{figure}

\begin{figure}[ht!]
     \includegraphics[width=0.47\textwidth]{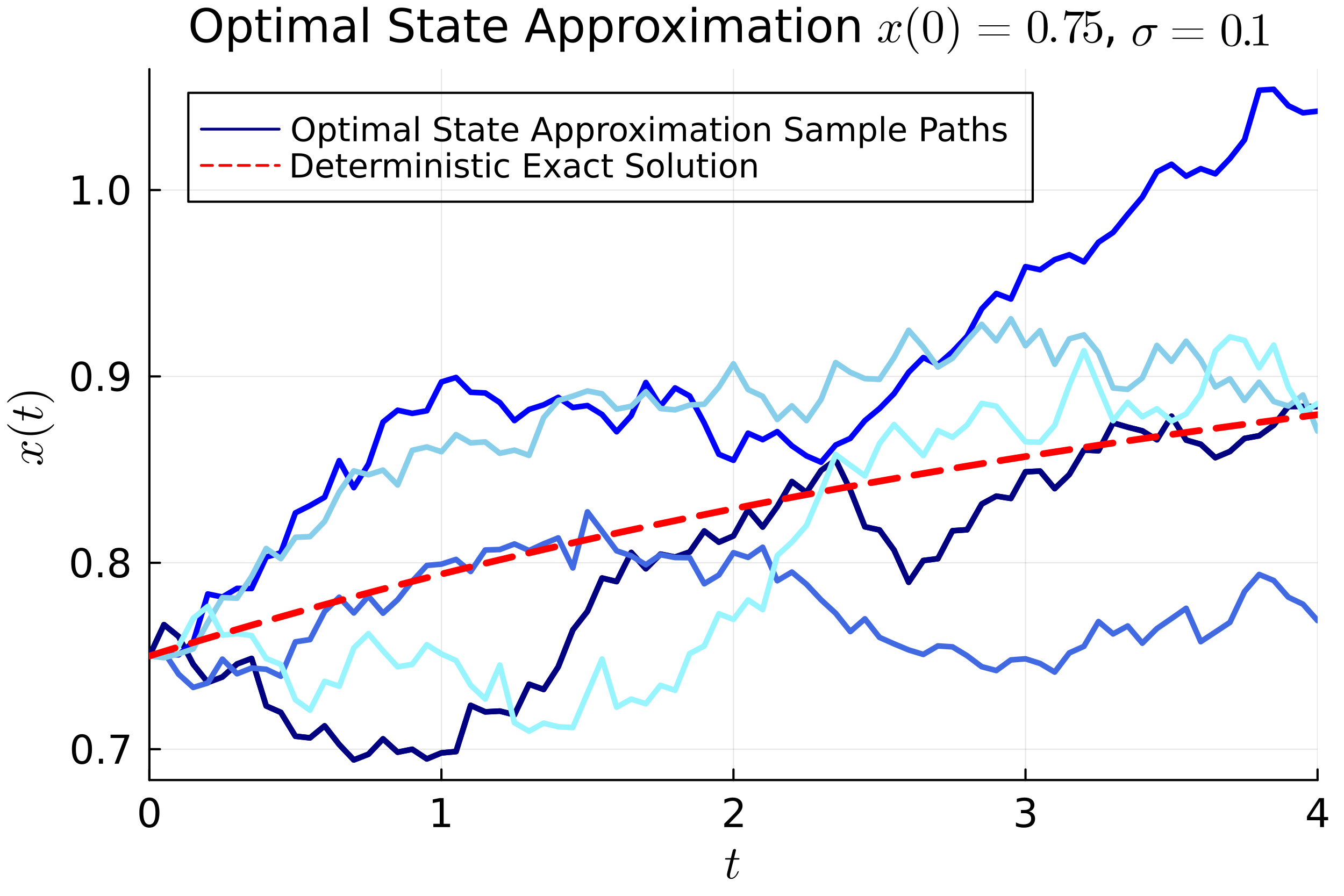}
     \qquad
     \includegraphics[width=0.47\textwidth]{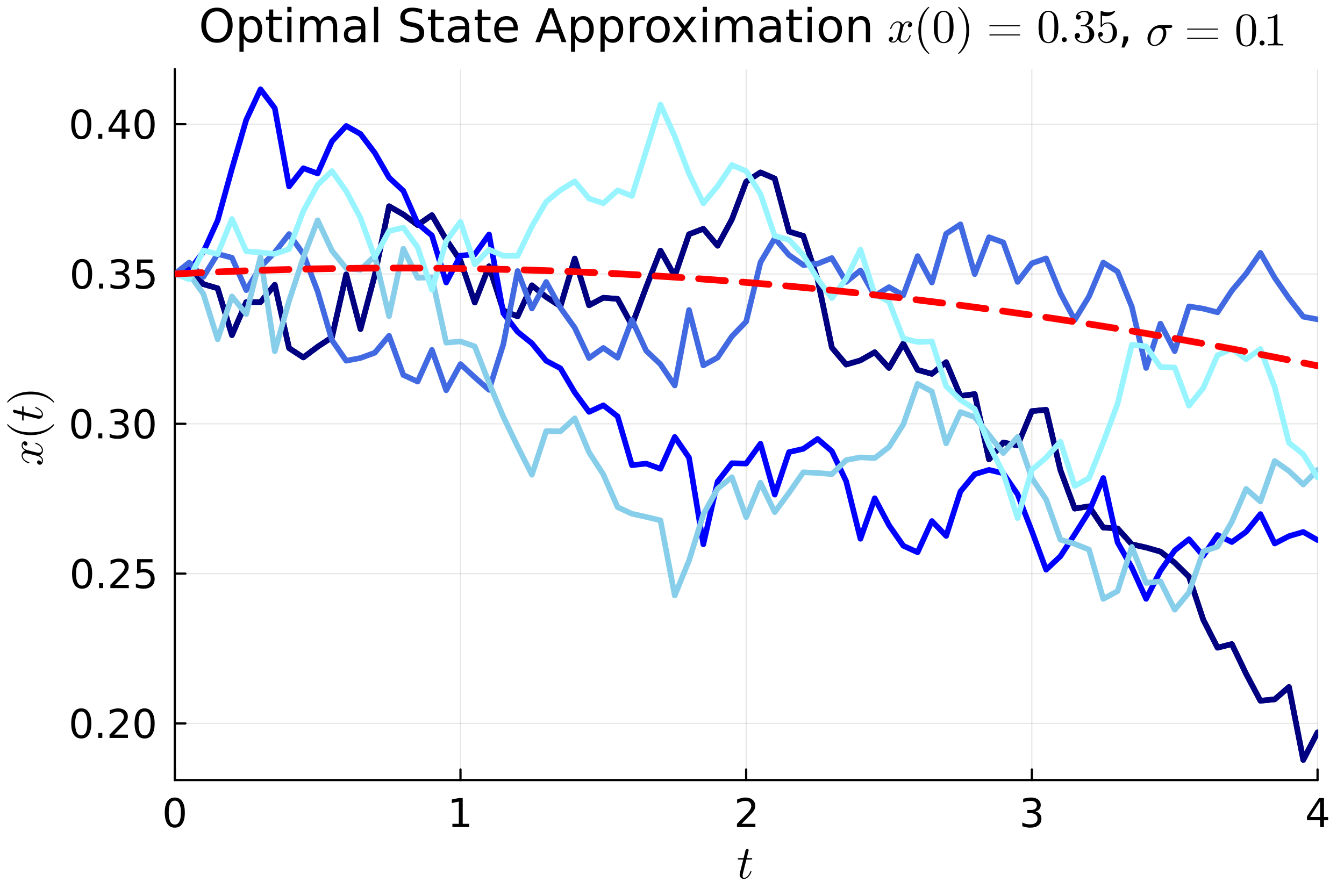}
     \caption{Deterministic optimal state trajectory and five sample paths of optimal state approximations with initial conditions $x(0)= 0.75\;(a)$ and $x(0)= 0.35\; (b)$, for $\sigma=0.1$.}
     \label{fig: State x0}
\end{figure}

In this section we present computational results for the two cases considered in~(\ref{subsec: Stochastic Dynamics}) and~(\ref{subsec: Stochastic State Observer}), obtained by solving the optimization problem in (\ref{eqn: NLP3}). For these experiments, we choose the following parameters values: $\alpha = 2$, $C = 1$, $\beta=\frac{1}{2}$, $\xi = \frac{1}{4}$, and $\sigma=0.1$. Figure~\ref{fig: Surfaces} shows, for both cases considered in~(\ref{subsec: Stochastic Dynamics}) and~(\ref{subsec: Stochastic State Observer}), respectively, results for two choices of number of Fourier coefficients: $M=N=3$ and $M=N=5$. To illustrate the behavior of sample approximating optimal trajectories around the deterministic one, Figure~\ref{fig: State x0} shows, for two fixed values of initial conditions, five samples of optimal state trajectory approximations obtained using the approximating optimal control.  

\section{Conclusions and Future Directions}\label{sec: Context and Future Directions}
This work presented a Fourier-based approximation for a family of optimal controls arising from a stochastic variation of the OCP presented in \cite{NFG22}. Our analysis and computational results showed that the control approximation is robust to the two variations of noise introduced in the dynamics of the system's state variable and that it drives the state trajectory up or down depending on the initial state $x(0)$, which means that past a certain proportion threshold, the optimal control, instead of driving the market-share up, slowly brings it down toward an equilibrium, which is otherwise not feasible in the absence of a control input. Future work will challenge some assumptions of the current method, such as the separability of state and control in the dynamics, which guarantees one unique set of Fourier coefficients in the parameterization of both control and state. The interpretability of the underlying descent in the space of Fourier coefficients will also be investigated, aiming for the potential extraction of theoretical properties of the intrinsic OCP from the coefficients convergence to the approximated optimal solution.   

\clearpage

\bibliographystyle{apsrev4-1}
\bibliography{References}

\begin{thebibliography}{10}%
\makeatletter
\providecommand \@ifxundefined [1]{%
 \@ifx{#1\undefined}
}%
\providecommand \@ifnum [1]{%
 \ifnum #1\expandafter \@firstoftwo
 \else \expandafter \@secondoftwo
 \fi
}%
\providecommand \@ifx [1]{%
 \ifx #1\expandafter \@firstoftwo
 \else \expandafter \@secondoftwo
 \fi
}%
\providecommand \natexlab [1]{#1}%
\providecommand \enquote  [1]{``#1''}%
\providecommand \bibnamefont  [1]{#1}%
\providecommand \bibfnamefont [1]{#1}%
\providecommand \citenamefont [1]{#1}%
\providecommand \href@noop [0]{\@secondoftwo}%
\providecommand \href [0]{\begingroup \@sanitize@url \@href}%
\providecommand \@href[1]{\@@startlink{#1}\@@href}%
\providecommand \@@href[1]{\endgroup#1\@@endlink}%
\providecommand \@sanitize@url [0]{\catcode `\\12\catcode `\$12\catcode
  `\&12\catcode `\#12\catcode `\^12\catcode `\_12\catcode `\%12\relax}%
\providecommand \@@startlink[1]{}%
\providecommand \@@endlink[0]{}%
\providecommand \url  [0]{\begingroup\@sanitize@url \@url }%
\providecommand \@url [1]{\endgroup\@href {#1}{\urlprefix }}%
\providecommand \urlprefix  [0]{URL }%
\providecommand \Eprint [0]{\href }%
\providecommand \doibase [0]{http://dx.doi.org/}%
\providecommand \selectlanguage [0]{\@gobble}%
\providecommand \bibinfo  [0]{\@secondoftwo}%
\providecommand \bibfield  [0]{\@secondoftwo}%
\providecommand \translation [1]{[#1]}%
\providecommand \BibitemOpen [0]{}%
\providecommand \bibitemStop [0]{}%
\providecommand \bibitemNoStop [0]{.\EOS\space}%
\providecommand \EOS [0]{\spacefactor3000\relax}%
\providecommand \BibitemShut  [1]{\csname bibitem#1\endcsname}%
\let\auto@bib@innerbib\@empty
\bibitem [{\citenamefont {Bass}(1969)}]{B69}%
  \BibitemOpen
  \bibfield  {author} {\bibinfo {author} {\bibfnamefont {F.~M.}\ \bibnamefont
  {Bass}},\ }\href@noop {} {\bibfield  {journal} {\bibinfo  {journal}
  {Management science}\ }\textbf {\bibinfo {volume} {15}},\ \bibinfo {pages}
  {215} (\bibinfo {year} {1969})}\BibitemShut {NoStop}%
\bibitem [{\citenamefont {Eryarsoy}\ \emph {et~al.}(2021)\citenamefont
  {Eryarsoy}, \citenamefont {Delen}, \citenamefont {Davazdahemami},\ and\
  \citenamefont {Topuz}}]{EDDT21}%
  \BibitemOpen
  \bibfield  {author} {\bibinfo {author} {\bibfnamefont {E.}~\bibnamefont
  {Eryarsoy}}, \bibinfo {author} {\bibfnamefont {D.}~\bibnamefont {Delen}},
  \bibinfo {author} {\bibfnamefont {B.}~\bibnamefont {Davazdahemami}}, \ and\
  \bibinfo {author} {\bibfnamefont {K.}~\bibnamefont {Topuz}},\ }\href@noop {}
  {\bibfield  {journal} {\bibinfo  {journal} {Journal of Business Research}\
  }\textbf {\bibinfo {volume} {124}},\ \bibinfo {pages} {163} (\bibinfo {year}
  {2021})}\BibitemShut {NoStop}%
\bibitem [{\citenamefont {Fan}\ \emph {et~al.}(2017)\citenamefont {Fan},
  \citenamefont {Che},\ and\ \citenamefont {Chen}}]{FCC17}%
  \BibitemOpen
  \bibfield  {author} {\bibinfo {author} {\bibfnamefont {Z.-P.}\ \bibnamefont
  {Fan}}, \bibinfo {author} {\bibfnamefont {Y.-J.}\ \bibnamefont {Che}}, \ and\
  \bibinfo {author} {\bibfnamefont {Z.-Y.}\ \bibnamefont {Chen}},\ }\href@noop
  {} {\bibfield  {journal} {\bibinfo  {journal} {Journal of business research}\
  }\textbf {\bibinfo {volume} {74}},\ \bibinfo {pages} {90} (\bibinfo {year}
  {2017})}\BibitemShut {NoStop}%
\bibitem [{\citenamefont {Kumar}\ and\ \citenamefont {Agarwala}(2016)}]{KA16}%
  \BibitemOpen
  \bibfield  {author} {\bibinfo {author} {\bibfnamefont {R.}~\bibnamefont
  {Kumar}}\ and\ \bibinfo {author} {\bibfnamefont {A.}~\bibnamefont
  {Agarwala}},\ }\href@noop {} {\bibfield  {journal} {\bibinfo  {journal}
  {Renewable and Sustainable Energy Reviews}\ }\textbf {\bibinfo {volume}
  {54}},\ \bibinfo {pages} {1515} (\bibinfo {year} {2016})}\BibitemShut
  {NoStop}%
\bibitem [{\citenamefont {Skiadas}\ and\ \citenamefont
  {Giovanis}(1997)}]{SG97}%
  \BibitemOpen
  \bibfield  {author} {\bibinfo {author} {\bibfnamefont {C.}~\bibnamefont
  {Skiadas}}\ and\ \bibinfo {author} {\bibfnamefont {A.}~\bibnamefont
  {Giovanis}},\ }\href@noop {} {\bibfield  {journal} {\bibinfo  {journal}
  {Applied stochastic models and data analysis}\ }\textbf {\bibinfo {volume}
  {13}},\ \bibinfo {pages} {85} (\bibinfo {year} {1997})}\BibitemShut {NoStop}%
\bibitem [{\citenamefont {Grasman}\ and\ \citenamefont
  {Kornelis}(2019)}]{GK19}%
  \BibitemOpen
  \bibfield  {author} {\bibinfo {author} {\bibfnamefont {J.}~\bibnamefont
  {Grasman}}\ and\ \bibinfo {author} {\bibfnamefont {M.}~\bibnamefont
  {Kornelis}},\ }\href@noop {} {\bibfield  {journal} {\bibinfo  {journal}
  {Journal of Mathematics in Industry}\ }\textbf {\bibinfo {volume} {9}},\
  \bibinfo {pages} {1} (\bibinfo {year} {2019})}\BibitemShut {NoStop}%
\bibitem [{\citenamefont {Fan}\ and\ \citenamefont {Griffin}(2017)}]{FG17}%
  \BibitemOpen
  \bibfield  {author} {\bibinfo {author} {\bibfnamefont {J.}~\bibnamefont
  {Fan}}\ and\ \bibinfo {author} {\bibfnamefont {C.}~\bibnamefont {Griffin}},\
  }\href@noop {} {\bibfield  {journal} {\bibinfo  {journal} {Operations
  Research Letters}\ }\textbf {\bibinfo {volume} {45}},\ \bibinfo {pages} {282}
  (\bibinfo {year} {2017})}\BibitemShut {NoStop}%
\bibitem [{\citenamefont {Nicolosi}\ \emph {et~al.}(2022)\citenamefont
  {Nicolosi}, \citenamefont {Friesz},\ and\ \citenamefont {Griffin}}]{NFG22}%
  \BibitemOpen
  \bibfield  {author} {\bibinfo {author} {\bibfnamefont {G.}~\bibnamefont
  {Nicolosi}}, \bibinfo {author} {\bibfnamefont {T.}~\bibnamefont {Friesz}}, \
  and\ \bibinfo {author} {\bibfnamefont {C.}~\bibnamefont {Griffin}},\
  }\href@noop {} {\bibfield  {journal} {\bibinfo  {journal} {Chaos, Solitons \&
  Fractals}\ }\textbf {\bibinfo {volume} {163}},\ \bibinfo {pages} {112535}
  (\bibinfo {year} {2022})}\BibitemShut {NoStop}%
\bibitem [{\citenamefont {Betts}(2020)}]{B20}%
  \BibitemOpen
  \bibfield  {author} {\bibinfo {author} {\bibfnamefont {J.~T.}\ \bibnamefont
  {Betts}},\ }\href@noop {} {\emph {\bibinfo {title} {Practical Methods for
  Optimal Control Using Nonlinear Programming}}},\ Vol.~\bibinfo {volume} {36}\
  (\bibinfo  {publisher} {SIAM},\ \bibinfo {year} {2020})\BibitemShut {NoStop}%
\bibitem [{\citenamefont {Schuster}\ and\ \citenamefont
  {Sigmund}(1983)}]{SS83}%
  \BibitemOpen
  \bibfield  {author} {\bibinfo {author} {\bibfnamefont {P.}~\bibnamefont
  {Schuster}}\ and\ \bibinfo {author} {\bibfnamefont {K.}~\bibnamefont
  {Sigmund}},\ }\href@noop {} {\bibfield  {journal} {\bibinfo  {journal}
  {Journal of theoretical biology}\ }\textbf {\bibinfo {volume} {100}},\
  \bibinfo {pages} {533} (\bibinfo {year} {1983})}\BibitemShut {NoStop}%
\end{thebibliography}%

\end{document}